\documentclass{amsart}
\usepackage{amsmath,amsfonts,amscd,amssymb}
\theoremstyle{plain}
\newcounter{thmcount}
\newtheorem{theorem}[thmcount]{Theorem}

\theoremstyle{definition}

\newtheorem{classification}[thmcount]{Classification}
\long\def\comment#1\endcomment{}
\input cyracc.def

\font\eightcyr=wncyr8
\def\smallsha{\text{\eightcyr\cyracc{Sh}}}

\def\F{{\mathbb F}}

\def\Q{{\mathbb Q}}
\def\R{{\mathbb R}}
\def\C{{\mathbb C}}

\def\P{{\mathbb P}}

\def\newmathop#1{\expandafter\gdef\csname #1\endcsname{\mathop{\rm #1}\nolimits}}
\newmathop{rk}
\newmathop{ord}
\newmathop{Gal}
\newmathop{SL}
\let\oldchar\char
\newmathop{char}

\let\char\oldchar
\newmathop{GL}
\newmathop{PGL}

\let\iso\cong
\let\tensor\otimes

\def\symk#1{\bigl(\frac{#1}k\bigr)}

\def\IZS{\text{\rm I$_0^*$}}

\def\II{\text{\rm II}}
\def\IIS{\text{\rm II$^*$}}
\def\III{\text{\rm III}}
\def\IIIS{\text{\rm III$^*$}}
\def\IV{\text{\rm IV}}
\def\IVS{\text{\rm IV$^*$}}

\begin{document}

\let\introdagger\dagger
\title{Elliptic curves with all quadratic twists of positive rank}
\author{Tim %$^\introdagger$
and Vladimir Dokchitser}
%\date{September 25, 2007}
%\subjclass[2000]{Primary 11G05; Secondary 11G07, 11G10, 11G40, 19A22, 20B99}
%\thanks{{\em MSC 2000:} Primary 11G05; Secondary 11G07, 11G40}
%\thanks{$^\dagger$Royal Society University Research Fellow}
\address{Robinson College, Cambridge CB3 9AN, United Kingdom}
\email{t.dokchitser@dpmms.cam.ac.uk}
\address{Gonville \& Caius College, Cambridge CB2 1TA, United Kingdom}
\email{v.dokchitser@dpmms.cam.ac.uk}
\llap{.\hskip 10cm} \vskip -8mm
\maketitle
\def\introdagger{{}}

\def\onehalf{{\raise3pt\hbox{$\scriptstyle 1$}\hskip-1.7pt{\hbox{\small/}}\hskip-1pt\raise-1pt\hbox{$\scriptstyle 2$}}}

%\vskip -9mm
%\llap{.\hskip 10cm}

%\section{Introduction}

Imagine you had an elliptic curve $E/K$ with everywhere good reduction,
defined over a number field $K$ that has no real and an odd number $r$
of complex places. Then the global root number $w(E/K)$ is $(-1)^r=-1$,
and it becomes $(-1)^{2r}=+1$ over every quadratic extension of $K$.
As the root number is the sign in the (conjectural) functional equation
for the $L$-function of $E$, the Birch--Swinnerton-Dyer conjecture
predicts that the Mordell-Weil rank of $E$ goes up in {\em every\/} quadratic
extension of $K$. Equivalently, every quadratic twist of $E/K$ has positive
rank, a behaviour that does not occur over $\Q$
(and would contradict Goldfeld's ``$\onehalf$ average rank'' conjecture).

These curves do exist\footnote{%
So Goldfeld's conjecture fails over number fields.
That there are curves all of whose quadratic twists must have positive
rank was observed in \cite{RohN} and is already implicit in \cite{Wal}.
% and, via Weil restriction,
%for abelian varieties over $\Q$.%
%So Goldfeld's conjecture fails for such $E/K$%
%and for their Weil restrictions to $\Q$.%
We also get, via Weil restriction, abelian varieties over $\Q$ all
of whose quadratic twists must have positive rank.
}.
For example, the elliptic curve over $\Q$
$$
% E: \quad y^2 + xy = x^3 + x^2 - 2x - 7 \qquad\qquad(\text{121C1})
  E: \quad y^2 = x^3 + \tfrac54 x^2 - 2x - 7 \qquad\qquad(\text{121C1})
$$
has discriminant $-11^4$ and acquires everywhere good reduction over any
cubic extension of $\Q$ which is totally ramified at 11. So
one may take
$K=\Q(\zeta_3,\sqrt[3]{11m})$ or $K=\Q(\sqrt[6]{-11m})$ for any positive
$m$, coprime to 11. (Those who prefer abelian extensions can take
$E=\text{1849C1}$ and $K$ to be the degree 6 field inside $\Q(\zeta_{43})$.)

%A surprising consequence for the curve $E=\text{121C1}$ is that
%for every $n\ge 1$ there is a number field $F$ where
%$$
%  \ord_{s=1} L(E_d/F,s) \ge n
%$$
%for every quadratic twist $E_d$ of $E/\Q$ by $d\in\Q^*$.
%!  This is NOT surprising: may take a multiquadratic extension of Q
%!  where every quad twist has (=> has lots) root numbers -1
%Indeed, take $F$ of the form
%$\Q(\zeta_3,\sqrt[3]{11m_1},\ldots,\sqrt[3]{11m_n})$;
%note also that here the $L$-functions are products of Rankin-Selberg
%convolutions and so are analytic at $s=1$.
%This contrasts the embarassingly unproven expectation
%that there is no number field $F$ where
%every elliptic curve over $\Q$ has rank $\ge 1$ (or $\ge 100$).

It is even easier to construct curves all of whose quadratic twists have
root number $+1$. For example,
$$
  E: \quad y^2 = x^3 + x^2 - 12x - \tfrac{67}4 \qquad\qquad(\text{1369E1})
$$
has discriminant $37^3$ and has everywhere good reduction over
$K=\Q(\sqrt[4]{-37})$. So it has root number $+1$ over every extension of $K$.
(Such a field $K$ exists for every elliptic curve with integral
$j$-invariant.) In view of the Birch--Swinnerton-Dyer conjecture, we expect
$E$ to have even rank
over {\em every\/} extension of $K$,
but it is not at all clear how to prove it for this or any other non-CM
elliptic curve%
\footnote{Similarly, as Karl Rubin remarked to us, there are fields $K$ such that $w(E/K)=1$
for {\em every\/}
elliptic curve $E$ defined over $\Q$; for instance $\Q(i,\sqrt{17})$ is such a field.}%
.

Let us say that an elliptic curve
$E/K$ is {\em lawful\/} if $w(E/K')=1$ for every quadratic extension $K'/K$,
and {\em chaotic\/} otherwise. Equivalently, $E$ is lawful if and only if
all of its quadratic twists have the same root number as $E/K$.
Depending on whether this root number is $+1$ or $-1$, let us call the curve
{\em lawful good\/} or {\em lawful evil}. Thus, conjecturally
the rank of a lawful evil curve increases in every quadratic extension%
\footnote{This is also implied by the conjectural finiteness of $\smallsha$,
under mild restrictions on $E$ at $v|6$, see \cite{Squarity} Thm. 1.3.}%
.

Another way of looking at our first example is that,
%\hbox{e.g.}
say,
for $K\!=\!\Q(\sqrt[6]{-11})$ the polynomial
$x^3\!+\!\tfrac54 x^2\!-\!2x\!-\!7$
must take {\em all\/} values%
\footnote{Can one prove (unconjecturally) that such square-free cubics exist,
over some $K$?
Note that there cannot be a non-constant parametric solution $(x(t),y(t))$ to $ty^2=f(x)$,
for otherwise $t\mapsto (x(t^2),t y(t^2))$ would be a non-constant map $\P^1\to E$.}
in $K^*/K^{*2}$.
Generally, one might conjecture that a square-free cubic $f(x)\!\in\! K[x]$
takes ``0\%'', ``50\%'' or all possible values in $K^*/K^{*2}$ depending
on whether the curve $y^2=f(x)$ is lawful good, chaotic or lawful evil over $K$.

\subsubsection*{Classification}

In our examples, we had the unnecessarily strong assumption that the
curve has everywhere good reduction.
Recall that the global root number $w(E/K)$ is the product of
local root numbers $w(E/K_v)$ over all places $v$ of $K$.
The condition that $w(E/K')=1$ for every quadratic
extension $K'/K$ is easily seen to be equivalent to $w(E/K'_v)$ being 1
for every quadratic extension $K'_v/K_v$, for all $v$.
For instance, this local condition fails for real places but holds for
complex places and primes of good reduction for $E$.

If $E$ is an elliptic curve over a local field $k$, let us also say that
$E/k$ is {\em lawful\/} if $w(E/k')=1$ for every quadratic extension $k'/k$,
and {\em chaotic\/} otherwise. Depending
on whether $w(E/k)$ is $+1$ or $-1$, call the curve
{\em lawful good\/} or {\em lawful evil}.
The reader should be warned that in the local setting lawfulness is not
equivalent to the invariance of the root number under quadratic twists.
%? Explain?

%If $K$ is a number field and $E/K$ is lawful evil,
%by the Birch--Swinnerton-Dyer conjecture
%the rank of $E$ should increase in every quadratic extension of $K$.

As mentioned above, a curve over a number field $K$ is lawful if
and only if it is
lawful over every completion of $K$. Whether it is good or evil
is determined by the parity of lawful evil places.
For instance, if $K$ has no real and $r$ complex places, an
elliptic curve $E/K$ with everywhere good reduction is lawful;
it is lawful evil if and only if $r$ is odd.

%\section{Lawfulness and locally abelian Galois action}

It turns out that if $E/K$ is lawful then $w(E/F)=w(E/K)^{[F:K]}$ for
{\em every\/} extension $F/K$.
%So the parity of the rank of $E/F$ should
%be determined by the parity of $[F:K]$ if $E/K$ is lawful evil.
In particular, a lawful evil curve $E/K$ must acquire points
of infinite order over any extension of even degree, while a lawful
good curve
should have even rank over every extension of $K$, as in the example
of 1369E1 above.
Generally,

\begin{theorem}
For an elliptic curve $E$ over a number field $K$,
the following conditions are equivalent:
\begin{enumerate}
\item $E/K$ is lawful, %i.e.{} $w(E/K')\!=\!1$ for every quadratic extension \hbox{$K'/K$,}
\item $E/K_v$ is lawful for all places $v$ of $K$,
\item $w(E/F)=w(E/K)^{[F:K]}$ for every finite extension $F/K$,
\item $K$ has no real places, and $E$ acquires everywhere good reduction
      over an abelian extension of $K$,
\item $K$ has no real places, and for all primes $p$ and all
places \hbox{$v\nmid p$} of~$K$,
the action of $\Gal(\bar K_v/K_v)$ on the Tate module $T_p(E)$ is abelian (``fake CM'')%
\footnote{Recall that $E$ has CM over $K$ if and only if the action of the {\em global\/}
Galois group $\Gal(\bar K/K)$ on $T_p(E)$ is abelian;
in this case it is obvious that $E$ has even rank over every extension of $K$.
}. %?
\end{enumerate}
\end{theorem}

\noindent
This is a corollary of the following local statement.

\begin{theorem}
For an elliptic curve $E$ over a local non-Archimedean field $k$
of characteristic 0, the following conditions are equivalent:
\begin{enumerate}
\item $E/k$ is lawful, %i.e.{} $w(E/k')\!=\!1$ for every quadratic extension \hbox{$k'/k$,}
\item $w(E/F)=w(E/k)^{[F:k]}$ for every finite extension $F/k$,
\item $E$ acquires good reduction over an abelian extension of $k$,
\item For some (any) $p$ different from the residue characteristic of $k$,
the action of $\Gal(\bar k/k)$ on $T_p(E)$ is abelian.
\end{enumerate}
%
%\noindent
%The implication
\end{theorem}

\begin{proof}
$(2)\Rightarrow(1)$ is obvious, and $(3)\Leftrightarrow(4)$
is a simple consequence of
the criterion of N\'eron--Ogg--Shafarevich in the case of
potentially good reduction and the theory of the Tate curve in the
potentially multiplicative case.
%
%\noindent
$(4)\Rightarrow(2)$ is an
elementary computation using the formula
$w(\chi)w(\bar\chi)=\chi(-1)$ (\cite{TatN} 3.4.7). %?
%? + norm vs restriction in local CFT
%(\cite{RohE}, p.145 or \cite{TatN} 3.4.7)
%
%\noindent
As for $(1)\Rightarrow(4)$,
if $k$ has odd residue characteristic or $E$ has potentially
multiplicative reduction, this follows from
the formulae for the local root numbers of Rohrlich~\cite{RohG} and
Kobayashi \cite{Kob}.
%? reformulate the above sentence?
In residue characteristic 2, computing root numbers is a nightmare.
Fortunately, Waldspurger has proved this in general on the other side
of the local Langlands correspondence:
a special or supercuspidal representation $\pi$ of $\PGL_2(k)$ has root
number $-1$ after base change to a suitable quadratic extension%
\footnote{If $E$ has potentially good
reduction and the Galois action on $T_p(E)$ is non-abelian,
the dual $V=(T_p(E)\tensor\Q_p)^*$ is an irreducible representation
of the Weil group of $k$. Twisting $V$ by the square root of the cyclotomic
character (this does not change the root number) gives a representation
that corresponds via local Langlands
to a supercuspidal representation $\pi$ of $\PGL_2(k)$.
%? Refer to J-L?
By \cite{Wal} Prop.~16, there is a quadratic extension $k'/k$ such
that $w(\pi/k')=-1$. Because the local Langlands correspondence
is known for $\PGL_2(k)$ (see \cite{JL}) and it preserves root numbers
and takes restriction to base change, we get $w(E/k')=-1$.}.
\end{proof}

An explicit classification of lawful elliptic curves in terms of
their $j$-invariant has been given by Connell, see \cite{Con} Prop. 6.
We end with an alternative classification in terms of the Kodaira symbols,
which also specifies whether the curve is good or evil.
It is easily deduced from the formulae for the local root numbers
of Rohrlich,
Kobayashi and the authors.
%, see \cite{RohG}~Thm. 2,~\cite{Kob}~Thm.~1.1 and \cite{Root2}.
%~\cite{RohG} Thm. 2,
%Kobayashi~\cite{Kob}~Thm.~1.1 and the authors~\cite{Root2}.

\begin{classification}
Let $k$ be a local field of characteristic 0, and $E/k$ an elliptic curve.
{}By \cite{RohG}~Thm. 2,

\begin{itemize}
\item
If $E$ has good reduction, then $E$ is lawful good.
\item
If $k\iso\C$, then $E$ is lawful evil.
\item
If $k\iso\R$ or $E$ has multiplicative or potentially multiplicative reduction,
  then $E$ is chaotic.
\end{itemize}

\noindent
Next, suppose $E/k$ has additive potentially good reduction and
$k$ has odd residue characteristic. Let $\Delta$ be the minimal
discriminant of $E$ of valuation $\delta$, and
write $\symk{\cdot}$ for the quadratic residue symbol on $k^*$.
{}From \cite{RohG}~Thm. 2 and \cite{Kob}~Thm.~1.1,
\begin{itemize}
\item
$E$ is lawful good
  if it has type $\IZS$ and $\symk{-1}=1$,\\
  type $\III, \IIIS$ with $\symk{-1}=1$ and $\symk{-2}=1$, or\\
  type $\II, \IIS, \IV, \IVS$ with $\Delta\in k^{*2}$ and
    $\symk{-1}^{\smash{\delta/2}}=1$.
\item
$E$ is lawful evil
  if it has type $\IZS$ and $\symk{-1}=-1$,\\
  type $\III, \IIIS$ with $\symk{-1}=1$ and $\symk{-2}=-1$, or\\
  type $\II, \IIS, \IV, \IVS$ with $\Delta\in k^{*2}$ and
    $\symk{-1}^{\smash{\delta/2}}=-1$.
\item $E$ is chaotic in all other cases.
\end{itemize}
(When $k$ has residue characteristic $\ge 5$, Kodaira types
\II, \III, \IV, \IZS, \IVS, \IIIS, \IIS{} correspond to $\delta=2,3,4,6,8,9,10$,
respectively. Also for $3\nmid\delta$, $\Delta\in k^{*2}$
%\newpage\noindent
if and only if
$\symk{-3}=1$ in this case, see \cite{Kob} 1.2.)

Finally, if $E/k$ has additive potentially good reduction and
$k$ has residue characteristic~2, let $c_4, c_6$ and $\Delta$
be the standard invariants
of some model of $E$, and set $\gamma(x)=x^8-6c_4x^4-8c_6x^2-3c_4^2\in k[x]$.
{}From \cite{Root2} Prop. 2 and Lemma 3,
$E/k$~is lawful if and only if
\begin{itemize}
\item $\sqrt{-3}\in k$ and $\gamma(x)$ is reducible, or
\item $\sqrt{-3}\not\in k$, $\sqrt[3]{\Delta}\in k$, and one of the
irreducible factors of $\gamma(x)$ becomes reducible over $k(\sqrt{-3})$.
\end{itemize}
In this case, $E$ is lawful good if and only if $-1$ is a norm from the
splitting field of $\gamma(x)$ to $k$ (\cite{Root2} Prop. 4b).
\end{classification}

%\medskip

Note from the classification that if $E/K$ has semistable reduction
at places above~2, it becomes lawful over some quadratic extension of $K$
if and only if it has integral $j$-invariant.
(If $v|2$ and $E/K_v$ has additive potentially good reduction,
$E$ stays chaotic in all quadratic extensions of $K_v$ if and only
the inertia group at~$v$ in $K(E[3])/K$ is $\SL_2(\F_3)$.)
In this way, it is easy to construct lawful evil curves over
imaginary quadratic fields. For example,
$$
  E: \quad y^2+xy=x^3-x^2-2x-1 \qquad\qquad(\text{49A1})
$$
is lawful evil over $\Q(i)$, so its rank should go up in every
extension of $\Q(i)$ of even degree.

\smallskip\noindent
{\bf Acknowledgements.}
We would like to thank David Rohrlich for suggesting to use
Waldspurger's results, and Andrew Granville for his interesting comments.
The first author is supported by a Royal Society
University Research Fellowship.

\end{document}